\documentclass[11pt]{article}
\usepackage{amsfonts,latexsym,rawfonts,amsmath,amssymb,amsthm,graphicx}
\numberwithin{equation}{section}
\newtheorem{theorem}{Theorem}[section]
\newtheorem{lem}[theorem]{Lemma}
\newtheorem{thm}[theorem]{Theorem}

\newtheorem{cor}[theorem]{Corollary}

\def\s{\,\,\,\,}

\def\endproof{$\hfill\Box$\\}
\def\R{\mathbb{R}}

\title{Addendum to "Contact stationary Legendrian surfaces in $\mathbb{S}^5$"[Pacific Math. J. 293(2018), no.1, 101-120]}
\author {Yong Luo }
\date{}
\begin{document}
\maketitle
\begin{abstract}
In \cite{Luo}, the present author proved that if $L$ is a contact stationary Legendrian surface in $\mathbb{S}^5$ with the canonical Sasakian structure and the square length of its second fundamental form belongs to $[0,2]$. Then we have that $L$ is either totally umbilical or is a flat minimal Legendrian torus. In this addendum we further prove that if $L$ is a totally umbilical contact stationary Legendrian surface in $\mathbb{S}^5$ , then $L$ is totally geodesic.
\end{abstract}

\section{Introduction}

In \cite{Luo},  we proved the following theorem:
\begin{thm}[\cite{Luo}]
Let $L:\Sigma\to \mathbb{S}^5$ be a contact stationary Legendrian surface. Then we have
\begin{eqnarray*}
\int_L\rho^2(3-\frac{3}{2}S+2H^2)d\mu\leq0,
\end{eqnarray*}
where $\rho^2:=S-2H^2$. In particular, if
\begin{eqnarray*}
0\leq S\leq 2,
\end{eqnarray*}
then either $\rho^2=0$ and $L$ is totally umbilical, or $\rho^2\neq 0$, $S=2, H=0$ and $L$ is a flat minimal Legendrian torus.
\end{thm}
Compared with the gap theorem of \cite{YKM}, it is very interesting to know if $L$ is totally geodesic in the above alternative when $\rho^2=0$. Hence
in the appendix of \cite{Luo}, we asked whether a totally umbilical contact stationary Legendrian surface in $\mathbb{S}^5$ with $0\leq S\leq 2$ is totally geodesic or not. In this note we give an affirmative positive answer to this question. Actually we get a stronger result.
\begin{thm}\label{main thm}
Assume that $L$ is a totally umbilical contact stationary Legendrian surface in $\mathbb{S}^5$. Then $L$ is totally geodesic.
\end{thm}
As a corollary of the above two theorems, we have
\begin{cor}
Assume that $L$ is a contact stationary Legendrian surface in $\mathbb{S}^5$ with $0\leq S\leq 2$. Then either $S=0$ and $L$ is totally geodesic or $S=2$ and $L$ is a flat minimal Legendrian torus.
\end{cor}
\section{Proof of Theorem \ref{main thm}}

Let $L$ be a Legendrian surface in $\mathbb{S}^5$ with the induced metric $g$. Assume that  $\{e_1,e_2\}$ is an orthonormal frame on $L$ such that $\{e_1,e_2,Je_1,Je_2,\textbf{R}\}$ be a orthonormal frame on $\mathbb{S}^5$. Here $\textbf{R}$ is the Reeb field of $\mathbb{S}^5$.

In the following we use indexes $i,j,k,l,s,t,m$ and $\beta,\gamma$ such that
\begin{eqnarray*}
1\leq i,j,k,l,s,t,m&\leq&2,
\\1\leq\beta,\gamma&\leq&3,
\\ \gamma^\ast=\gamma+2,\s \beta^\ast&=&\beta+2.
\end{eqnarray*}

Let $B$ be the second fundamental form of $L$ in  $\mathbb{S}^5$ and define
\begin{eqnarray}
h_{ij}^k&=&g_\alpha(B(e_i,e_j),Je_k),
\\h^3_{ij}&=&g_\alpha(B(e_i,e_j),\textbf{R}).
\end{eqnarray}
Then
\begin{eqnarray}
h_{ij}^k&=&h_{ik}^j=h_{kj}^i,
\\h^3_{ij}&=&0.
\end{eqnarray}
The Gauss equations and Ricci equations are
\begin{eqnarray}
R_{ijkl}&=&(\delta_{ik}\delta_{jl}-\delta_{il}\delta_{jk})+\sum_s(h^s_{ik}h^s_{jl}-h^s_{il}h^s_{jk}),\label{basic equation 1}
\\R_{ik}&=&\delta_{ik}+2\sum_sH^sh^s_{ik}-\sum_{s,j}h^s_{ij}h^s_{jk},
\\2K&=&2+4H^2-S,
\\R_{3412}&=&\sum_i(h_{i1}^1h_{i2}^2-h_{i2}^1h_{i1}^2)\nonumber
\\&=&\det h^1+\det h^2,
\end{eqnarray}
where $K$ is the sectional curvature function of $(L,g)$ and $h^1,h^2$ are the second fundamental forms w.r.t. the normal directions $Je_1$, $Je_2$ respectively.

 In addition we have the following Codazzi equations and Ricci identities
\begin{eqnarray}
h^\beta_{ijk}&=&h^\beta_{ikj},
\\h^\beta_{ijkl}-h^\beta_{ijlk}&=&\sum_mh^\beta_{mj}R_{mikl}+\sum_mh^\beta_{mi}R_{mjkl}+\sum_\gamma h^\gamma_{ij}R_{\gamma^\ast\beta^\ast kl}.\label{basic equation 2}
\end{eqnarray}

Using these equations, we can get the following Simons' type inequality:
\begin{lem}[\cite{Luo}]\label{main result}
Let $L$ be a Legendrian surface in $\mathbb{S}^5$. Then we have
\begin{eqnarray}\label{main lemma}
\frac{1}{2}\Delta\sum_{i,j,\beta}(h^\beta_{ij})^2&\geq&|\nabla^T h|^2-2|\nabla^T H|^2-2|\nabla^\nu H|^2 +\sum_{i,j,k,\beta}(h^\beta_{ij}h^\beta_{kki})_j \nonumber
\\&+&S-2H^2+2(1+H^2)\rho^2-\rho^4-\frac{1}{2}S^2,
\end{eqnarray}
where $|\nabla^T h|^2=\sum_{i,j,k,s}(h^s_{ijk})^2$ and $|\nabla^T H|^2=\sum_{i,s}(H^s_i)^2$.
\end{lem}
\proof This lemma was proved in \cite{Luo}. We copy the proof here because we will use several equalities and inequalities in the proof in the following. Using equations from (\ref{basic equation 1}) to (\ref{basic equation 2}), we have
\begin{eqnarray}\label{simon type}
\frac{1}{2}\Delta\sum_{i,j,\beta}(h^\beta_{ij})^2
&=&\sum_{i,j,k,\beta}(h^\beta_{ijk})^2+\sum_{i,j,k,\beta}h^\beta_{ij}h^\beta_{kijk}\nonumber
\\&=&|\nabla h|^2-4|\nabla^\nu H|^2+\sum_{i,j,k,\beta}(h^\beta_{ij}h^\beta_{kki})_j+\sum_{i,j,l,k,\beta} h^\beta_{ij}(h^\beta_{lk}R_{lijk}+h^\beta_{il}R_{lj})\nonumber
\\&+&\sum_{i,j,k,\beta,\gamma} h^\beta_{ij}h^\gamma_{ki}R_{\gamma^\ast\beta^\ast jk}\nonumber
\\&=&|\nabla h|^2-4|\nabla^\nu H|^2+\sum_{i,j,k,s}(h^s_{ij}h^s_{kki})_j+2K\rho^2-2(\det h^1+\det h^2)^2\nonumber
\\&\geq&|\nabla h|^2-4|\nabla^\nu H|^2+\sum_{i,j,k,\beta}(h^\beta_{ij}h^\beta_{kki})_j+2(1+H^2)\rho^2-\rho^4-\frac{1}{2}S^2,
\end{eqnarray}
where $\rho^2:=S-2H^2$ and in the above calculations we used the following identities
\begin{eqnarray*}
\sum_{i,j,k,l,\beta} h^\beta_{ij}(h^\beta_{lk}R_{lijk}+h^\beta_{il}R_{lj})&=&2K\rho^2,
\\\sum_{i,j,k,\beta,\gamma} h^\beta_{ij}h^\gamma_{ki}R_{\gamma^\ast\beta^\ast jk}&=&-2(\det h^1+\det h^2)^2,
\end{eqnarray*}
where in the first equality we used $R_{lijk}=K(\delta_{lj}\delta_{ik}-\delta_{lk}\delta_{ij})$ and $R_{lj}=K\delta_{lj}$ in a proper orthonormal frame field, because $L$ is a surface.

Note that
\begin{eqnarray}\label{main idea1}
|\nabla h|^2&=&\sum_{i,j,k,\beta}(h^\beta_{ijk})^2\nonumber
=|\nabla^T h|^2+\sum_{i,j,k}(h^3_{ijk})^2\nonumber
=|\nabla^T h|^2+\sum_{i,j,k}(h^k_{ij})^2\nonumber
\\&=&|\nabla^T h|^2+S,
\end{eqnarray}
where in the third equality we used
\begin{eqnarray*}
h^3_{ijk}&=&\langle\bar{\nabla}_{e_k}h(e_i,e_j),\textbf{R}\rangle
\\&=&-\langle h(e_i,e_j),\bar{\nabla}_{e_k}\textbf{R}\rangle
\\&=&\langle h(e_i,e_j),Je_k\rangle
\\&=&h^k_{ij}.
\end{eqnarray*}
Similarly we have
\begin{eqnarray}\label{main idea2}
|\nabla^\nu H|^2=|\nabla^TH|^2+H^2.
\end{eqnarray}
Combing (\ref{simon type}), (\ref{main idea1}) with (\ref{main idea2}), we get (\ref{main lemma}).
\endproof

We also have
\begin{lem}[\cite{Luo}]
Let $L:\Sigma\to \mathbb{S}^5$ be a contact stationary Legendrian surface. Then
\begin{eqnarray}\label{integral equality}
\int_L|\nabla^\nu H|^2d\mu=-\int_L(K-1)H^2d\mu,
\end{eqnarray}
where $|\nabla^\nu H|^2=\sum_{\beta,i}(H^\beta_i)^2$.
\end{lem}

Integrating over (\ref{main lemma}) and using $|\nabla^Th|^2\geq 3|\nabla^TH|^2$ (see appendix, Lemma A.1 of \cite{Luo}) we get
\begin{eqnarray}\label{ine1}
0&\geq&\int_L[(|\nabla^T h|^2-3|\nabla^T H|^2)-2|\nabla^\nu H|^2+S-2H^2+2(1+H^2)\rho^2-\rho^4-\frac{1}{2}S^2+|\nabla^T H|^2]d\mu \nonumber
\\ &\geq& \int_L[-2|\nabla^\nu H|^2+S-2H^2+2(1+H^2)\rho^2-\rho^4-\frac{1}{2}S^2+|\nabla^T H|^2]d\mu \nonumber
\\&=&\int_L(2-\rho^2)\rho^2d\mu+\int_L 2H^2\rho^2+2(K-1)H^2-2H^2+S-\frac{1}{2}S^2+|\nabla^T H|^2d\mu \nonumber
\\&=&\int_L(2-\rho^2)\rho^2d\mu+\int_L 2H^2\rho^2+(4H^2-S)H^2-2H^2+S-\frac{1}{2}S^2+|\nabla^T H|^2d\mu \nonumber
\\&=&\int_L\frac{3}{2}\rho^2(2-S)+2H^2\rho^2+|\nabla^T H|^2d\mu,
\end{eqnarray}
where in the first equality we used (\ref{integral equality}) and in the second equality we used the Gauss equation $2K=2+4H^2-S$.

Therefore we obtain the following integral inequality
\begin{eqnarray}\label{ine2}
\int_L\rho^2(3-\frac{3}{2}S+2H^2)+|\nabla^T H|^2d\mu\leq0.
\end{eqnarray}
Particularly if  $\rho^2=0$, i.e.  $L$ is totally umbilical, then from (\ref{ine2}) we see that $|\nabla^T H|^2=0$. Then from (\ref{main idea2}) we get that $|\nabla^\nu H|^2=H^2$, which implies that $\int_LKH^2d\mu=0$ by (\ref{integral equality}). Now by the Gauss equation $2K=2+4H^2-S=2+2H^2-\rho^2=2+2H^2$ we get
$$\int_LH^2(1+H^2)=0.$$
Therefore $H=0$ and hence combing with the assumption that $0=\rho^2=S-2H^2$,  we get $S=0$, i.e. $L$ is totally geodesic.

This completes the proof of Theorem \ref{main thm}.
\endproof

\vspace{1cm}

\textbf{Acknowledgement.} The author is supported by the NSF of China(No.11501421).
{}
\vspace{1cm}\sc
Yong Luo

School of Mathematics and statistics, Wuhan University, Wuhan 430072, China.

{\tt yongluo@whu.edu.cn}

\vspace{1cm}\sc
\end{document}